# About Factorial Sums

Mihály Bencze[1] and Florentin Smarandache[2]
[1]Str. Hărmanului 6, 505600 Săcele-Négyfalu, Jud. Braşov, Romania
[2]Chair of Math & sciences, University of New-Mexico, 200 College Road, NM 87301, USA

**Abstract.** In this paper, we present some new inequalities for factorial sum.

**Application 1.** We have the following inequality
$$\sum_{k=1}^{n} k! \le \frac{2((n+1)!-1)}{n+1}$$

**Proof.** If $x_k$, $y_k > 0$ $(k=1,2,...,n)$, have the same monotonity, then
$$\left(\frac{1}{n}\sum_{k=1}^{n} x_k\right)\left(\frac{1}{n}\sum_{k=1}^{n} y_k\right) \le \frac{1}{n}\sum_{k=1}^{n} x_k y_k \qquad (1)$$
the Chebishev's inequality.

If $x_k$, $y_k$ have different monotonity, then holds true the reverse inequality, we take $x_k = k$, $y_k = k!$ $(k=1,2,...,n)$ and use that $\sum_{k=1}^{n} k \cdot k! = (n+1)! - 1$.

**Application 2.** We have the following inequality
$$\sum_{k=1}^{n} k! \le \frac{3(n+1)(n+1)!}{n^2 + 3n + 5}$$
**Proof.** In (1) we take
$$x_k = k^2 + k + 1 ;$$
$$y_k = k! \ (k=1,2,...,n)$$
and the identity
$$\sum_{k=1}^{n} (k^2 + k + 1) k! = (n+1)(n+1)!$$

**Application 3.** We have the following inequality
$$\sum_{k=1}^{n} \frac{1}{k!} \ge \frac{n^2(n+1)}{2((n+1)!-1)}$$
**Proof.** Using the Application 1, we take

$$\sum_{k=1}^{n}\frac{1}{k!} \geq \frac{n^2}{\sum_{k=1}^{n}k!} \geq \frac{n^2(n+1)}{2((n+1)!-1)}$$

**Application 4.** We have the following inequality
$$\sum_{k=1}^{n}\frac{1}{k!} \geq \frac{n^2(n^2+3n+5)}{3(n+1)(n+1)!}$$

**Proof.** Using the Application 2, we take
$$\sum_{k=1}^{n}\frac{1}{k!} \geq \frac{n^2}{\sum_{k=1}^{n}k!} \geq \frac{n^2(n^2+3n+5)}{3(n+1)(n+1)!}$$

**Application 5.** We have the following inequality:
$$\sum_{k=1}^{n}\frac{1}{k!} \geq 1+\frac{2}{n}\left(1-\frac{1}{n!}\right)$$

**Proof.** In (1) we take $x_k = k$, $y_k = \frac{1}{(k+1)!}$, $(k=1,2,...,n)$ and we obtain

$$\frac{1}{n}\left(\sum_{k=1}^{n}k\right)\left(\sum_{k=1}^{n}\frac{1}{(k+1)!}\right) \geq \sum_{k=1}^{n}\frac{k}{(k+1)!} = 1-\frac{1}{(n+1)!}$$

therefore
$$\left(\sum_{k=1}^{n}\frac{1}{(k+1)!}\right) \geq \frac{2}{n+1}\left(1-\frac{1}{(n+1)!}\right)$$

or
$$\sum_{k=2}^{n}\frac{1}{k!} \geq \frac{2}{n}\left(1-\frac{1}{n!}\right)$$

therefore
$$\left(\sum_{k=1}^{n}\frac{1}{k!}\right) \geq 1+\frac{2}{n}\left(1-\frac{1}{n!}\right)$$

**Application 6.** We have the following inequality:
$$\left(\sum_{k=1}^{n}\frac{1}{(k+2)^2 k!}\right) \geq \frac{2}{n+5}\left(1-\frac{1}{(n+2)!}\right)$$

**Proof.** In (1) we take $x_k = k+2$, $y_k = \frac{1}{(k+2)^2 k!}$, $(k=1,2,...,n)$

therefore
$$\frac{1}{n}\left(\sum_{k=1}^{n}(k+2)\right)\sum_{k=1}^{n}\frac{1}{(k+2)^2 k!} \geq \sum_{k=1}^{n}\frac{1}{(k+2)^2 k!} = 1-\frac{1}{(n+2)!}$$

therefore

$$\sum_{k=1}^{n} \frac{1}{(k+2)^2 k!} \geq \frac{2}{n+5}\left(1 - \frac{1}{(n+2)!}\right)$$

**Application 7.** We have the following inequality:
$$\sum_{k=1}^{n} \frac{1}{k(k+1)(k+2)!} \geq \frac{6}{2n^2 + 9n + 1}\left(\frac{1}{2} - \frac{1}{(n+1)(n+2)!}\right)$$

**Proof.** In (1) we take
$$x_k = k^2 + 2k + 2, \quad y_k = \frac{1}{k(k+1)(k+2)!}, \quad (k = 1, 2, \ldots, n)$$

then

$$\frac{1}{n}\sum_{k=1}^{n}(k^2 + 2k + 2)\sum_{k=1}^{n}\frac{1}{k(k+1)(k+2)!} \geq \sum_{k=1}^{n}\frac{k^2 + 2k + 2}{k(k+1)(k+2)!} =$$
$$= \sum_{k=1}^{n}\frac{1}{k(k+1)!} - \frac{1}{(k+1)(k+2)!} = \frac{1}{2} - \frac{1}{(n+1)(n+2)!}$$

**Application 8.** We have the following inequality:
$$\sum_{k=1}^{n} \frac{1}{4k^4 + 1} \geq \frac{n}{2n^2 + 2n + 1}$$

**Proof.** In (1) we take $x_k = 4k$, $y_k = \frac{1}{4k^4 + 1}$, $(k = 1, 2, \ldots, n)$,

therefore
$$\frac{1}{n}\left(\sum_{k=1}^{n} 4k\right)\left(\sum_{k=1}^{n} \frac{1}{4k^4 + 1}\right) \geq \sum_{k=1}^{n} \frac{4k}{4k^4 + 1} = \sum_{k=1}^{n}\left(\frac{1}{2k^2 - 2k + 1} - \frac{1}{2k^2 + 2k + 1}\right) = \frac{2n(n+1)}{2n^2 + 2n + 1}$$

**Application 9.** We have the following inequality:
$$\sum_{k=1}^{n} \frac{1}{4k^4 - 1} \geq \frac{3n}{(2n+1)^2}$$

**Proof.** In (1) we take $x_k = k^2$, $y_k = \frac{1}{4k^2 - 1}$, $(k = 1, 2, \ldots, n)$ then

$$\frac{1}{n}\left(\sum_{k=1}^{n} k^2\right)\left(\sum_{k=1}^{n} \frac{1}{4k^2 - 1}\right) \geq \sum_{k=1}^{n} \frac{k^2}{4k^2 - 1} = \frac{n(n+1)}{2(2n+1)}, \text{ etc.}$$

**Reference:**

[1] Octogon Mathematical Magazine (1993-2007)

{Published in Octogon Mathematical Magazine, Vol. 15, No. 2, 810-812, 2007.}